\documentclass[12pt]{article}
\usepackage{amssymb}
\usepackage{amsthm}
\usepackage{amsmath}
\usepackage{tikz}
\usetikzlibrary{calc}
\DeclareSymbolFont{extraup}{U}{zavm}{m}{n}
\DeclareMathSymbol{\clubsuit}{\mathalpha}{extraup}{88}

\newtheorem{thm}{Theorem}[section]
\newtheorem{lem}[thm]{Lemma}
\newtheorem{prop}[thm]{Proposition}
\newtheorem{conj}{Conjecture}

\newtheorem{prob}[conj]{Problem}

\title{Finding one excellent element in case of one lie}
\author{{\bf Aanchal Gupta}\\Delhi, India\\aang1995@gmail.com
\and {\bf Gyula O.H. Katona}\thanks{The work of the second author was
supported by the Hungarian National Research, Development and Innovation Office
NKFIH under grant numbers SSN135643 and K132696. }\\
R\'enyi Institute, HUN-REN\\
Budapest Pf 127, 1364 Hungary \\ ohkatona@renyi.hu  }

\begin{document}
\date{}

\maketitle

\begin{abstract} An $n$-element set contains an unknown number of "excellent" elements and we want to find one. The members of a family ${\cal F}\subset 2^{[n]} $ of subsets can be asked if they contain at least one excellent element or not. At most one of the answers can be wrong. We find the smallest family which finds one excellent element or claims that there is none.


\end{abstract}

\section{Introduction}

Suppose that we have $n$ rock samples, some of them are radioactive, and we need
one radioactive sample. A possible test is taking an arbitrary subset and
checking if this set is radiating, that is, if there is at least one radioactive
among them. The mathematical problem, of course, is the minimum number
of tests sufficient to find at least one of the radioactive ones.

Let $[n] = \{ 1, 2, . . . , n\}$ be an $n$-element set. Some of them have certain special properties. In the traditional model of Combinatorial Search Theory this property is usually negative, they are mostly called {\it defective}, following the motivating applications. In these examples usually we need to find all defective elements using the questions  of type ”does the subset $F \subset [n] $
contain at least one defective element?”. The following survey papers and
books contain many pretty examples, applications and results of this type:
\cite{K1}, \cite{A}, \cite{AW}, \cite{DH}, \cite{GP}.

Ray Chambers \cite{Ch} raised the question what happens if not all defective
elements have to be found, only one and the number of defective elements is unknown. 
In this situation the word {\it excellent} is more proper for the special property.
This will be used in the rest of the paper. The problem is denoted by $P_n(?\rightarrow 1)$.

In the case of every search problem there are two essentially different
approaches. If the search is adaptive then the choice of the subsequent questions may depend on the answers to the previous questions. More precisely,
first a set $F$ is tested. If the answer is Y, that is $F$ contains at least one
excellent element then a set $F_1$ is tested afterwards, while in the case when
the answer is N, that is $F$ contains no excellent element then another set $F_0$
is tested. This is repeated. If e.g. the answer to the question $F_1$ is Y then
the next question is $F_{11}$, and so on. This is continued until the goal of the
search is achieved, in the present paper when at least one excellent element
is determined, or it is established that there is none. The question sets form
a tree-like structure. This is called a search algorithm. Of course the number
of questions may depend on the actual placement of the excellent elements.
The number of questions in the worst case (length of the longest path in
the tree) is the length of the search. The {\it adaptive} complexity of the search
problem is the minimum of the length of the search taking for all search
algorithms solving the given problem.
The non-adaptive search, however, is just a family ${\cal F}$ of the question sets.
They are all tested and the goal of the search should be achieved based on
these answers. The length of the {\it non-adaptive} search is just $|{\cal F}|$. The non-adaptive complexity of the search problem is $\min |{\cal F}|$ for all families solving
the given problem. Since the non-adaptive search is also an adaptive search,
the adaptive complexity cannot exceed the non-adaptive complexity for any
search problem.

The following easy statement can be found in the literature cited above.

\begin{prop} If it is known that there is exactly one excellent element then both the adaptive and non-adaptive complexities are $\lceil \log_2 n \rceil$.
		
\end{prop}

Paper \cite{K2} proved the following statements.

\begin{thm}
The adaptive complexity of the problem $P_n(?\rightarrow 1)$ is  $\lceil \log_2 (n+1) \rceil$,
while its non-adaptive complexity is $n$.
\end{thm}

A combinatorial search {\it algorithm of two rounds}
 consists of a family ${\cal F} \subset 2^{[n]}, |{\cal F}|=m$
and a set of families ${\cal G}(s) \subset 2^{[n]}$ defined for every 
$s \in \{ N, Y\}^m$.
The members
of ${\cal F}$  are asked in the first round. If the sequence of answers is $s$ then the
members of ${\cal G}(s)$ are asked in the second round.
 The length of this algorithm is
$|{\cal F}| + max_s |{\cal G} (s)|$.

\begin{thm} {\rm \cite{K2}} The two-round complexity of the the problem $P_n(?\rightarrow 1)$
	is $\sim 2\sqrt{n}$.

\end{thm}

Gerbner and Vizer \cite{GV} generalized this for $r$ rounds: $\sim rn^{1\over r}$.

\section{The result} 

Until this point it was supposed that the results of the tests, the answers to the questions are correct. R\'enyi \cite{R} and Ulam \cite{U} independently suggested to study the combinatorial search problems in the case when some of the answers might be incorrect.
This is called {\it R\'enyi-Ulam game} or {\it search in presence of lies}. For a good survey see \cite{D}. Here we will suppose that there is at most one lie. In the case of non-adaptive search it means the following. The algorithm is a family ${\cal F}\subset 2^{[n]}$ of subsets, we ask if $F\in {\cal F}$ contains an excellent element or not, for every member of the family. All answers are correct with at most one exception. But, unlike in the previous case, a set $F$ can be repeated, so ${\cal F}$ is actually  a {\it multi-family}, where the members have multiplicities 1, 2 or 3. The algorithm is good if either an excellent element was found on the basis of these answers or we can claim that there is none. Later in the proofs we will speak about the {\it disjoint partition} of a multi-family. If the set$F$  occurs twice 
(three times) the the ese copies might be separated: one of them goes to obe subfamily the other one(s) goes to the other subfamily.

\begin{thm} The non-adaptive complexity of the problem $P(?\rightarrow 1)$ in case of at most one lie  is $2n+1$.
\end{thm}

{\it Proof of the upper bound: the optimal algorithm}

Let ${\cal F}$ consists of $[n]$ and all sets $\{ i\}$ with multiplicity 2. Let us show that we find a excellent element for all possible sequences of answers.

{\it Case 1}. There is an $\{ i\}$ such that we obtain two different answers for the two copies.

One of them is {\it the} lie, all other answers are correct. That is the two answers are either YES, YES, or NO, NO for all other elements.

{\it Case 1.1}. The answer for the question $[n]$ is YES.

If we got two answers YES for the two copies of $\{ j\}$ then $j$ is a excellent element.
If the answers for all $\{ j\} (\not= i)$ are NO, NO, then $i$ is the excellent element.

{\it Case 1.2}. The answer for the question $[n]$ is NO.

This answer must be correct, therefore there is not excellent element in this case.

{\it Case 2}. For every $i$ the two answers for the question $\{ i\}$ are the same.

{\it Case 2.1}. There is an $i$ with answers YES, YES.

Then independently of the answer for $[n]$, $i$ is a excellent element.

{\it Case 2.2}. The two answers are NO, NO for all $\{ i\}$.

Then independently of the answer for $[n]$, there is no excellent element.
\vskip 5mm
{\it Proof of the lower bound}

\begin{lem} Suppose that the multi-family ${\cal F}\subset 2^{[n]}$ solves the problem 
	$P(?\rightarrow 1)$ in the case where at most one answer can be wrong
	if and only if the following two conditions hold.
	
	(i) Every element is covered by at least 3 members of  ${\cal F}$.
	
	(ii) Take a disjoint partition of ${\cal F}$ of into two subfamilies ${\cal F}_1$ and ${\cal F}_2$ where $|{\cal F}_1|\geq 2$, define
	$$T^{*}({\cal F}_1, {\cal F}_2)=
	\{ T:\ T\cap F\not= \emptyset , T\cap G=\emptyset$$  $${\mbox{ holds for every }} F\in {\cal F}_1, G\in {\cal F}_2 {\mbox{ with at most one exception}}\} , $$
	suppose that this family is not empty then
	$$\bigcap_{T\in T^{*}({\cal F}_1, {\cal F}_2)}T\not=\emptyset.\eqno(1)$$
\end{lem}

{\it Proof}

{\it Necessity of (i).} Let $x\in [n]$, and suppose that $x\in F\in {\cal F}, x\in G\in {\cal F}$ and $x$ is not an element of the other members of ${\cal F}$. If the answers are NO for all questions but $F$ (for which the answer is YES) then we cannot decide if there is no excellent element at all (the answer for $F$ was a lie) or $x$ may be excellent (the answer for $G$ was a lie).

{\it Necessity of (ii).} The problem is solved for any sequence of answers. Suppose that the answers for the questions $F\in {\cal F}_1$ are YES and the answers for the questions $G\in {\cal F}_2$ are NO. Denote the set of  excellent elements by $T$. Then $T\cap F\not= \emptyset , T\cap G=\emptyset$ for every $F\in {\cal F}_1$ and $G\in {\cal F}_2$ must hold with at most one exception. In other words $T$ is a member of the family $T^{*}({\cal F}_1, {\cal F}_2)$. Since the answer is YES for at least two sets and one must be a correct answer, $T\not=\emptyset$ can be supposed. We can be sure that $x$ is excellent only if it is an element of every member of this family. This implies that the intersection of all family members is non-empty.

{\it Sufficiency.} Suppose that (i) and (ii) hold for ${\cal F}$. If the answer for YES is obtained for none or only one member of ${\cal F}$ then by (i) for any $x$ there are at least two
members saying NO, consequently there is no excellent element. On the other hand if there are at least two members of ${\cal F}$ with answer YES then let this subfamily be ${\cal F}_1$.
By (ii) there is an $x$ which is excellent. $\Box$

\begin{lem} (ii) implies that ${\cal F}$ has a one-element member.
\end{lem}	
	
	{\it Proof}. Suppose the contrary: all members contain at least two elements. Apply (ii) with ${\cal F}_1={\cal F}$. Then
	$$[n]-\{ i\} \in T^{*}({\cal F}, \emptyset )$$
	holds for every $i\in [n]$ since this set intersects every subset except $\{ i\}$ which is not in ${\cal F}$, by our assumption.   Then  $\cap_{i=1}^n([n]-\{ i\} )=\emptyset$ contradicts (1).
 $\Box$

\begin{lem} (ii) implies that there is a $j\in [n]$ such that $\{ j\} $ is a member of
	${\cal F}$ with multiplicity $u\geq 2$.
\end{lem}

{\it Proof.} Suppose that there is no such $j$.  By the previous lemma we have a 
$\{ j\}\in {\cal F}$, but it has multiplicity one. Let $i\in [n]$ be an arbitrary element. Use (ii), again, with
${\cal F}_1={\cal F}$. Then
$$[n]-\{ i\} \in T^{*}({\cal F}, \emptyset )\eqno(2)$$
holds when $[n]-\{ i\}$ intersects every member of ${\cal F}$ with one exception.
If the set $\{ j\}\in {\cal F}$ has multiplicity one, then taking this as an exception,
$[n]-\{ i\}$ meets every member of ${\cal F}-\{ j\}$, therefore (2) holds for every
$i\in [n]$. The intersection of these sets is empty, contradicting (1). $\Box$
\vskip 5mm
{\it End of the proof of the lower estimate}

We prove that $|{\cal F}|\geq 2n+1$ holds for every non-adaptive algorithm on $n$ elements by induction on $n$. The case $n=1$ is trivial. Suppose that the statement is true for $n-1$ and prove it for $n$.

By Lemma 2.4 we know that there is a set $\{ j\}$ with multiplicity $u\geq 2$ in  ${\cal F}$. For sake of simplicity suppose that $j=n$. Define the families
$${\cal F}^{\prime}=\{ F:  F\in {\cal F}, n\not\in F \}.$$ 
$${\cal F}^-=\{ F-\{n \}: F\in {\cal F}, n\in F,  F\not=\{ n\} \},$$
$${\cal F}[n-1]={\cal F}^{\prime}\cup {\cal F}^-.$$
Note here that if $F\in {\cal F}^{\prime}, F\in {\cal F}^-$ both hold then the multiplicities of $F$ add up in ${\cal F}[n-1]$.

Since we deleted at least two members (the multiplicity $u$ of $\{ n\}$  in ${\cal F}$ can be more than 2) we have
$$|{\cal F}|\geq |{\cal F}[n-1]|+2. \eqno(3)$$
Let us show that ${\cal F}[n-1]$ satisfies (i) and (ii) for the underlying set $[n-1]$.

Let $i\not=n$. We know by Lemma 2.2 that ${\cal F}$ satisfies (i) on the underlying set $[n]$ therefore there are 3 members of ${\cal F}: F_1, F_2, F_3$ such that $i\in F_1, F_2, F_3$.
Then   $i\in F_1-\{ n\}, F_2-\{n \}, F_3-\{ n\}$. But these are members of ${\cal F}[n-1]$, (i) holds.

Let us prove that (ii) holds for ${\cal F}[n-1]$ on the underlying set $[n-1]$. 

Let ${\cal G}_1\cup {\cal G}_2$ be a decomposition of ${\cal F}[n-1]$ where $|{\cal G}_1|\geq 2$. Define
$${\cal G}_i^{\prime}=\{ G\in {\cal G}_i\cap {\cal F}^{\prime}\}$$
and
$${\cal G}_i^{-}=\{ G\cup \{n \}: G\in {\cal G}_i\cap {\cal F}^{-}\}$$
for $i=1,2$. It is clear that ${\cal G}_1^{\prime}, {\cal G}_1^{-}, {\cal G}_2^{\prime}$ and ${\cal G}_2^{-}$ are pairwise disjoint subfamilies of ${\cal F}$. Moreover their union is ``almost" ${\cal F}$, namely only the copies of $\{ n\}$ are missing. Consequently
$$({\cal G}_1^{\prime}\cup {\cal G}_1^{-})\bigcup ({\cal G}_2^{\prime}\cup {\cal G}_2^{-}\cup ({\mbox{ copies of}} \{ n\}))$$
is a decomposition of ${\cal F}$ where $|{\cal G}_1^{\prime}\cup {\cal G}_1^{-}|\geq 2$.
By (ii) we know that the intersection of all sets $T$ satisfying
$$T\in T^{*}({\cal G}_1^{\prime}\cup {\cal G}_1^{-},{\cal G}_2^{\prime}\cup {\cal G}_2^{-}\cup ({\mbox{ copies of}} \{ n\} ))\eqno(4)$$
is non-empty. Since $T$ must be disjoint to all but one member of ${\cal G}_2^{\prime}\cup {\cal G}_2^{-}\cup {\mbox{ copies of} }\ \{ n\}$ and the multiplicity of $\{ n\}$ is at least 2,
we obtain $n\not\in T$. Therefore the conditions in the definition of $T^{*}$ can be restricted to $[n-1]$. If the members of ${\cal F}$ are restricted to $[n-1]$ then
${\cal G}_i^{\prime}\cup {\cal G}_i^{-}$ becomes ${\cal G}_i$, (4) is equivalent to
$$T\in T^{*}({\cal G}_1, {\cal G}_2),$$
therefore the intersection of all sets $T$ is non-empty, as desired, (ii) is satisfied for
${\cal F}[n-1]$. By Lemma 2.2  ${\cal F}[n-1]$ is a non-adaptive algorithm solving the prroblem $P_{n-1}(?\rightarrow 1)$ in the case of at most one lie.  By the inductional hypothesis $|{\cal F}[n-1]|\geq 2(n-1)+1$. Using (3) we obtain $|{\cal F}|\geq 2n+1.$
$\Box$

\section{Conclusions}
An easy non-adaptive algorithm shows that $2n+1$ questions are sufficient to find one excellent element if their number is unknown and there can be one lie among the answers. It might be surprising that it is only the double of the value in the case when there is no lie.  We proved that one cannot go below $2n+1$. 

But there are many related unsolved open question here.
\begin{prob} Determine the non-adaptive complexity of the problem $P_n(?\rightarrow 1)$
in the case of at most $t$ lies.
\end{prob}
\begin{prob} Determine the adaptive complexity of the problem $P_n(?\rightarrow 1)$ in the case of at most 1 lie.
\end{prob}
\begin{prob} What happens if we want to find more excellent elements, say $s$ of them or to claim that there are only $r<s$. The answer is not known even in the case when no lie is allowed.
\end{prob}

\end{document}